\documentclass[11pt,reqno,oneside]{amsart}

\usepackage[utf8]{inputenc}
\usepackage[T1]{fontenc}

\usepackage{fouriernc}

\usepackage{graphics}
\usepackage[pagebackref]{hyperref}
\renewcommand*\backref[1]{\ifx#1\relax \else (Cited on page #1) \fi}

\usepackage{amsmath,amsthm,amssymb}

\usepackage{graphics}
\usepackage{hyperref}
\usepackage[dvipsnames]{xcolor}

\usepackage{booktabs}

\usepackage{mathtools}
\mathtoolsset{showonlyrefs}

\usepackage{enumitem}
\newlist{todolist}{itemize}{2}
\setlist[todolist]{label=$\square$}
\usepackage{pifont}
\usepackage[square,sort,comma,numbers]{natbib}

\definecolor{darkblue}{rgb}{0.0,0.0,0.3}
\hypersetup{colorlinks,breaklinks,
  linkcolor=darkblue,urlcolor=darkblue,
anchorcolor=darkblue,citecolor=darkblue}

\theoremstyle{plain}
\newtheorem{theorem}{Theorem}
\newtheorem*{theorem*}{Theorem}

\newtheorem*{proposition*}{Proposition}

\newtheorem*{corollary*}{Corollary}

\theoremstyle{definition}
\newtheorem{remark}[theorem]{Remark}

\newtheorem*{algorithm*}{Algorithm}

\newtheorem*{experiment*}{Experiment}

\DeclareMathOperator{\SL}{SL}
\DeclareMathOperator{\OLD}{OLD}

\DeclareMathOperator*{\EE}{\mathbb{E}}

\newcommand{\todaysdate}{\the\year.\ifnum\month<10 0\fi\the\month.\ifnum\day<10 0\fi\the\day}

\title{On Murmurations and Trace Formulas}
\author{David Lowry-Duda \\ 30 May 2025}
\date{}

\begin{document}

\begin{abstract}
  In recent work with Bober, Booker, Lee, Seymour-Howell, and Zubrilina, we
  proved murmuration behavior for Maass forms in the eigenvalue aspect and for
  modular forms in the weight aspect.
  Both used an approach based on the Selberg trace formula.
  But different trace formulas, including those due to Kuznetsov or Petersson,
  offer different variations.
  We examine murmurations from the perspective of different trace formulas and
  outline several families of $L$-functions where one can likely prove
  additional murmuration behavior.
\end{abstract}

\maketitle

\section{Introduction}

Using output from machine learning algorithms, He, Lee, Oliver, and
Pozdnyakov~\cite{murmur2022} noticed that root numbers of $L$-functions
associated to elliptic curves seem to correlate with the traces of Frobenius
$a_p$ at primes $p$ in proportion to the conductor.
They dubbed this a ``murmuration''.
The language of murmurations is $L$-functions.
A general (self-dual) $L$-function has the shape
\begin{equation}
  L(s, \pi) = \sum_{n \geq 1} \frac{\lambda(n)}{n^s}
\end{equation}
and satisfies a functional equation
\begin{equation}\label{eq:feq}
  \Lambda(s, \pi)
  :=
  N^s G(s) L(s, \pi)
  =
  \varepsilon
  N^{1-s} G(1-s) L(1-s, \pi),
\end{equation}
where $N$ is the conductor of $L(s, \pi)$, $G(s)$ is a product of $\Gamma$
functions and powers of $\pi$, and $\varepsilon = \pm 1$ is the root number.
(See \S5.1 of~\cite{iwanieckowalski04} for a more complete, but still very
general description).

Murmuration phenomena are correlations between the coefficients $\lambda(n)$
and root numbers $\varepsilon$ for $L$-functions associated to $\pi$ in some
family $\mathcal{F}$.
In practice, murmuration phenomena seem to be the strongest when the family
$\mathcal{F}$ is ordered by conductor (though as Sutherland noted in his
talk~\cite{sutherland24SCGPmurmuration}, different choices of ordering
can lead to weaker---but present!---murmurations).

\section*{Acknowledgements}

This note expands on previous talks of the
author~\cite{dld23heuristicmaassmurmur, dld24scgpmurmur}
from the Hot Topics
\href{https://icerm.brown.edu/program/hot_topics_workshop/htw-23-ma}{Murmurations
in Arithmetic} workshop at ICERM 6--8 July 2023 and the
\href{https://scgp.stonybrook.edu/archives/43116}{Murmurations in Arithmetic
Geometry and Related Topics} workshop at the Simons Center for Geometry and
Physics 11--15 November 2024.
I'm grateful to the organizers of both workshops.
I thank Jonathan Bober, Andy Booker, Alex Cowan, Min Lee, Peter Sarnak, Andrei
Seymour-Howell, and Nina Zubrilina for helping make murmurations fun and
enriching.

This work was supported by the Simons Collaboration in Arithmetic
Geometry, Number Theory, and Computation via the Simons Foundation grant 546235.

\section{General Murmurations and $L$-functions}\label{sec:general_murmurations}

To describe general murmurations, it is useful to use the notational framework
from~\cite{sarnak23letter}.
Given a smooth nonnegative weight function $\Phi: (0, \infty) \longrightarrow
\mathbb{R}$ of compact support, a family $\mathcal{F}$ of $L$-functions $L(s,
\pi)$ ordered by conductor $N(\pi)$, and a function $f: \mathcal{F}
\longrightarrow \mathbb{C}$, we define
\begin{equation}
  A_\Phi^\mathcal{F}(f, X)
  =
  A(f, X)
  :=
  \sum_{\pi \in \mathcal{F}}
  \Phi\big( N(\pi) / X \big) f(\pi),
\end{equation}
as well as the expected value
\begin{equation}
  \EE_{\pi \in \mathcal{F}}[f; X; \Phi]
  =
  \EE_{\pi \in \mathcal{F}}
  [f; X]
  :=
  \frac{A(f, X)}{A(1, X)}.
\end{equation}
In this notation, we recover the initial murmurations of~\cite{murmur2022} by
taking $\Phi$ to be the indicator function on $[1, 2]$, $\mathcal{F}$ to be the
family $\mathcal{E}^{\pm}$ of elliptic curves with root number $\varepsilon =
\pm 1$, and $f$ to be $a_E(p)$, the trace of Frobenius at $p$.
This gives
\begin{equation}
  \EE_{\mathcal{E}^{\pm}} [a_E(p); X]
  =
  \frac{%
    \displaystyle \sum_{\substack{E \in \mathcal{E}^{\pm} \\ X \leq N(E) \leq 2X}} a_E(p)
  }{%
    \displaystyle \sum_{\substack{E \in \mathcal{E}^{\pm} \\ X \leq N(E) \leq 2X}} 1
  }.
\end{equation}
The plots in~\cite{murmur2022} suggest correlations between the conductors and
$\mathcal{E}_{\mathcal{E}^\pm}[a_E(p); X]$ for various primes $p$ (or really,
for various ratios $p/N(E)$).

This notation allows a more precise definition of murmuration phenomena.
If there exists a function $M_\Phi: (0, \infty) \longrightarrow \mathbb{R}$
such that
\begin{equation}\label{eq:murmuration_shape}
  \EE_{\pi \in \mathcal{F}}[\lambda_\pi(p); X; \Phi]
  =
  M_\Phi(p/X)
  +
  \mathrm{Err}_\Phi(p, N)
\end{equation}
where $\mathrm{Err}_\Phi(p, N)$ is a smaller error term, then $M_\Phi$
describes a statistical correlation between the coefficients $\lambda_\pi(p)$
and the conductors $N(\pi)$.
We call this a \emph{murmuration phenomenon}.

\begin{remark}[On normalization]
  The relationship between normalization and murmurations can be overstated.
  For example, in~\cite{BBLLD23} the murmuration density is given for both
  $\lambda_\pi(p)$ and $\lambda_\pi(p)\sqrt{p}$; they're similar, except
  that multiplying by $\sqrt{p}$ boosts the magnitude of the murmuration
  pattern.

  For elliptic curves over $\mathbb{Q}$, the traces of Frobenius $a_E(p)$ and
  the coefficients of the analytically-normalized $L$-function $\lambda_E(p)$
  are related via $a_E(p) = \lambda_E(p) \sqrt{p}$.
  On average, these satisfy $\lvert a_E(p) \rvert^2 \asymp p$ and $\lvert
  \lambda_E(p) \rvert^2 \asymp 1$.
  As the signs of $a_E(p)$ vary approximately randomly, one should expect the
  size of $A(a_E(p), X)$ to behave like a random walk with steps of size
  $\sqrt{p}$.
  Such a walk experiences square-root cancellation, so that one should expect
  $A(a_E(p), X) \asymp \# \{ E : X \leq N(E) \leq 2X \}$, matching the size of
  $A(1; X)$.
  Thus square-root normalization leads to murmuration patterns of constant
  size.

  But omitting the scaling factor (or by undernormalizing or overnormalizing)
  doesn't necessarily remove the statistical correlation.
\end{remark}

\section{One-Level Density}

In his letter~\cite{sarnak23letter}, Sarnak notes that~\cite{ils2000lowlying}
includes a proof of an arithmetically weighted murmuration for the family of
self-dual $\SL(2, \mathbb{Z})$ modular forms of increasing weight.
We describe this in some detail as it motivates further work.

Katz and Sarnak~\cite{katzsarnak99zeroes, katz1999random} predict that
the statistics of zeros in certain families of $L$-functions are
governed by certain random matrix models from groups of matrices associated to
the symmetry types of the $L$-functions.
One such set of statistics is the one-level density, which concerns the
distribution of low-lying zeros of $L$-functions.

We continue the notation from \S\ref{sec:general_murmurations}.
Let $\mathcal{F}(N) := \{ \pi \in \mathcal{F} : N(\pi) = N \}$ denote
elements of the family $\mathcal{F}$ with conductor $N$ and let $\phi(x)$
be an even function of rapid decay as $\lvert x \rvert \to \infty$.
The one-level density of $\mathcal{F}$ is
\begin{equation}
  \OLD_\phi(\mathcal{F})
  :=
  \lim_{N \to \infty} \frac{1}{\# \mathcal{F}(N)}
  \sum_{\pi \in \mathcal{F}(N)}
  \sum_{\gamma_\pi} \phi \Big( \frac{\gamma_\pi \log N}{2 \pi} \Big),
\end{equation}
where $\gamma_\pi$ runs through zeros of
$\Lambda(\frac{1}{2} + i \gamma_\pi, \pi)$, counted with multiplicity.
The normalizing factor of $\gamma_\pi$ guarantees that the zeros have
unit spacing on average.
As $\phi$ decays rapidly, the one-level density describes the density of zeros
of $L(s, \pi)$ within $O(1 / \log N)$ of the critical point with varying
$\phi$ measures.

According to the Katz-Sarnak philosophy, for many families there is a measure
$W_\mathcal{F}$ coming from matrices such that
\begin{equation}
  \OLD_\phi(\mathcal{F})
  =
  \int_{\mathbb{R}}
  \widehat{\phi}(x)
  \widehat{W_\mathcal{F}}(x)
  dx
\end{equation}
for all nice test functions $\phi$.
The first main theorem of~\cite{ils2000lowlying} is one such result for weight
functions $\phi$ of limited support.

\begin{theorem*}[Theorem 1.1 of~\cite{ils2000lowlying}]
  Assume GRH.\@
  Let $\phi$ be an even Schwarz function with
  $\mathrm{supp}(\widehat{\phi}) \subset (-2, 2)$.
  Let $H_k^{\pm}(N)$ denote a Hecke eigenbasis of modular newforms of level
  $1$, weight $k$, and root number $\varepsilon = \pm 1$.
  Then
  \begin{equation}\label{eq:ilsOLD}
    \OLD_\phi(H_k^{\pm})
    =
    \int_\mathbb{R}
    \widehat{\phi}(x)
    \widehat{W_{\mathrm{SO}(\pm)}}(x)
    dx
  \end{equation}
  where
  $W_{\mathrm{SO}(+)}(x)
  =
  W_{\mathrm{SO}(\mathrm{even})}(x)
  = 1 + \frac{\sin(2 \pi x)}{2 \pi x}$
  and
  $W_{\mathrm{SO}(-)}(x)
  =
  W_{\mathrm{SO}(\mathrm{odd})}(x)
  =
  1 - \frac{\sin(2 \pi x)}{2 \pi x} + \delta_0(x)$,
  where $\delta_0(x)$ is the Kronecker $\delta$ symbol at $0$.
\end{theorem*}
The distributions $W_{\mathrm{SO}(\mathrm{odd})}$ and
$W_{\mathrm{SO}(\mathrm{even})}$, as well as distributions corresponding to
certain other matrix groups, were computed in~\cite[Theorem~AD.12.6]{katz1999random}.

Unravelling, this theorem shows (under GRH) that
\begin{equation}\label{eq:ils_expanded}
  \lim_{N \to \infty}
  \frac{1}{\# H_k^{\pm}(N)}
  \sum_{f \in H_k^{\pm}(N)}
  \sum_{\gamma_f}
  \phi \Big(
    \frac{\gamma_f \log N}{2 \pi}
  \Big)
  =
  \int_\mathbb{R}
  \widehat{\phi}(x)
  \widehat{W_{\mathrm{SO}}}(\pm)(x)
  dx
\end{equation}
The explicit formula relating $L$-function zeros to sums over primes
relates the sums
\begin{equation}\label{eq:explicit_formula}
  \sum_{\gamma_f}
  \phi \Big(
    \frac{\gamma_f \log N}{2 \pi}
  \Big)
  \leftrightsquigarrow
  \sum_p
  \frac{\lambda_f(p) \log p}{\sqrt p}
  \widehat{\phi}
  \Big( \frac{\log p}{\log N} \Big).
\end{equation}
This use is approximately as in~\cite[\S5.5]{iwanieckowalski04}; see \S4
of~\cite{ils2000lowlying} for the precise relationship.

Note that if $\widehat{\phi}$ is supported on $[-\theta, \theta]$, then only
primes $\leq N^\theta$ contribute.
Thus inserting~\eqref{eq:explicit_formula} into~\eqref{eq:ils_expanded} and
swapping the order of summation shows that
$\OLD_\phi(H_k^{\pm})$ behaves like
\begin{equation}\label{eq:double_exp}
  \EE_{p \sim N^\theta}
  \EE_{f \in H_k^{\pm}(N)}
  [\lambda_f(p) \log p / \sqrt{p}].
\end{equation}
This almost describes murmuration behavior with some prime scaling, except that
there is averaging around primes $p \sim N^\theta$.
Murmurations concern the behavior around primes $p \sim N$.

Remarkably, the behavior of murmurations occurs exactly at the transition range
of the one-level density.
More specifically, the Fourier transforms of
$W_{\mathrm{SO}(\mathrm{odd})}$ and $W_{\mathrm{SO}(\mathrm{even})}$ satisfy
\begin{align}
  \widehat{W_{\mathrm{SO}(\mathrm{odd})}}(y) &=
  \delta_0(y) + \frac{\mathbf{1}_{[-1, 1]}(y)}{2},
  \\
  \widehat{W_{\mathrm{SO}(\mathrm{even})}}(y) &=
  \delta_0(y) + \frac{2 - \mathbf{1}_{[-1, 1]}(y)}{2}.
\end{align}
There is a discontinuity in behavior in the one-level density~\eqref{eq:ilsOLD}
exactly when $\pm 1 \in \mathrm{supp}(\widehat{\phi})$, which occurs exactly
when primes relevant to murmurations begin to contribute.

Although~\cite{ils2000lowlying} doesn't directly state a murmuration result,
the proofs imply an arithmetically weighted murmuration.
Specifically, let $L(s, \mathrm{Sym}^2 f)$ denote the symmetric square
$L$-function associated to the modular form $f$ and let $S^{\pm}_k(1)$ denote a
basis of Hecke cuspforms with root number $\pm 1$ for the space of weight $k$,
level $1$ modular forms.
(Actually the root numbers are determined by $k \bmod 4$. If $k \equiv 0 \bmod
4$, then the root number is $1$. If $k \equiv 2 \bmod 4$, then the root number
is $-1$. There are no level $1$ cuspforms if $k \not \equiv 0 \bmod 2$.)
Note that $L(1, \mathrm{Sym}^2 f)$ is positive~\cite[\S5.12]{iwanieckowalski04}.
Then the proofs in~\cite{ils2000lowlying} imply (still assuming GRH) that
\begin{equation}\label{eq:ils_murmuration}
  \frac{%
    \sum_{f \in S_k^{\pm}} \Phi(N(f) / X)
    \frac{\lambda_f(p) \sqrt{p}}{L(1, \mathrm{Sym}^2 f)}
  }{%
    \sum_{f \in S_k^{\pm}} \Phi(N(f) / X)
    \frac{1}{L(1, \mathrm{Sym}^2 f)}
  }
  =
  \pm 4 \pi
  \sum_{c \geq 1}
  \frac{\mu^2(c)}{c^2 \varphi(c)}
  \Phi \left(
    \frac{16 \pi^2 p}{c^2 X}
  \right)
  +
  \mathrm{Err}_\Phi(p, X),
\end{equation}
in which $\mu(\cdot)$ is the Möbius function, $\varphi(\cdot)$ is Euler's
totient function, and the error term $\mathrm{Err}_\phi(p, X)$ is dominated by
the main term as $X \to \infty$.
The only difference between~\eqref{eq:ils_murmuration} and a murmuration of the
standard shape~\eqref{eq:murmuration_shape} is the arithmetic weights
$L(1, \mathrm{Sym}^2 f)^{-1}$.
Unfortunately, it's not apparent how to translate an arithmetically weighted
murmuration to a non-arithmetically weighted normalization.

Prompted by a suggestion of Sarnak,
in~\cite{BBLLD23} Booker, Bober, Lee, and the author sought to study
non-arithmetically normalized murmurations of $\mathrm{SL}(2, \mathbb{Z})$
modular cuspforms of increasing weight.
After performing an arbitrarily small average over primes, one can show a
murmuration of the following form.
\begin{theorem*}[Theorem~1.1 of~\cite{BBLLD23}]
  Assume GRH.\@
  Fix a compact interval $E \subset \mathbb{R}_{> 0}$ with
  $\lvert E \rvert > 0$.
  Set $N = \mathcal{N}(K) \approx (\tfrac{k-1}{4\pi})^2$, the conductor of a
  modular cuspform of weight $K$.
  As $K \to \infty$,
  \begin{equation}
    \frac{
      \sum_{\substack{p\;\textup{prime} \\ p/N \in E}}
      \log p
      \sum_{k \sim K}
      \sum_{f \in S_k^{\pm}}
      \lambda_f(p)
    }{
      \sum_{\substack{p\;\textup{prime} \\ p/N \in E}}
      \log p
      \sum_{k \sim K}
      \sum_{f \in S_k^{\pm}}
      1
    }
    =
    \frac{\pm 1}{\sqrt{N}}
    \Big(
      \frac{\nu(E)}{\lvert E \rvert}
      +
      o_{E, \epsilon}(1)
    \Big),
  \end{equation}
  where
  \begin{equation}\label{eq:nu}
    \nu(E)
    =
    \frac{1}{\zeta(s)}
    \sideset{}{^\ast}\sum_{
      \substack{
        a, q \in \mathbb{Z}_{> 0}
        \\ \gcd(a, q) = 1
        \\ (a / q)^{-2} \in E
      }
    }
    \frac{\mu(q)^2}{\varphi(q)^2 \sigma(q)}
    \big(
      \frac{q}{a}
    \big)^3,
  \end{equation}
  where the $\ast$ indicates that terms occurring at endpoints of $E$ should be
  halved.
\end{theorem*}

The additional small average over the interval $E$ makes the notation
cumbersome, but it is remarkable how similar the murmuration functions are
in~\eqref{eq:ils_murmuration} and~\eqref{eq:nu}; notably, there are point
masses at squares of squarefree integers (visible in~\eqref{eq:nu} after taking
$a = 1$).

\section{Trace formulas}

Although the distributions are different, the path to the results
in~\cite{ils2000lowlying} and~\cite{BBLLD23} share a similar structure: they
both start with a trace formula.
Iwaniec, Luo, and Sarnak use the Petersson trace
formula~\cite{petersson32trace}:
\begin{equation}\label{eq:petersson}
  \Big(
    \frac{\Gamma(k-1)}{(4 \pi \sqrt{mn})^{k-1}}
  \Big)
  \sum_{f \in S_k}
  \frac{\lambda_f(n) \lambda_f(m)}{\lVert f \rVert^2}
  =
  \delta_{[m = n]}
  +
  2 \pi i^k
  \sum_{c > 0}
  \frac{S(m, n; c)}{c}
  J_{k-1} \Big(
    \frac{4 \pi \sqrt{mn}}{c}
  \Big),
\end{equation}
where $\lVert f \rVert^2 = \langle f, f \rangle$ is the Petersson inner
product, $S(m, n; c)$ is the Kloosterman sum
\begin{equation}
  \sum_{d \in (\mathbb{Z} / c \mathbb{Z})^\times}
  \exp \Big(
    \frac{2 \pi i (md + n d^{-1})}{c}
  \Big),
\end{equation}
and $J_\nu(x)$ is the Bessel function of the first kind.
To study murmurations, one takes $m = 1$, $n = p$, and studies the right-hand
side.
The arithmetic weights $L(1, \mathrm{Sym}^2 f)$ in~\eqref{eq:ils_murmuration}
come from the denominators $\lVert f \rVert^2$ in~\eqref{eq:petersson}.
As the coefficients $\lambda_f(n)$ are real for cuspforms of level $1$, the
special value $L(1, \mathrm{Sym}^2 f)$ is proportional to $\lVert f
\rVert^2$~\cite[eq.\ 5.101]{iwanieckowalski04}.

In contrast,~\cite{BBLLD23} applies a variant of the
Eichler--Selberg  
trace formula from~\cite{child22}, which gives an expression for the trace of
the $n$th Hecke operator on all of $S_k$.
(The downside of
Eichler--Selberg  
is that quadratic class numbers appear, which
are more mysterious than the Kloosterman sums from the Petersson trace
formula).

More generally, almost every established murmuration result currently
known to the author starts with a trace formula.
We briefly describe them.

\begin{enumerate}
  \item Lee, Oliver, and Pozdnyakov used Gauss sums and Poisson summation to
  study murmurations of Dirichlet characters in~\cite{murmurDirichlet2023}.
  (Sarnak's letter~\cite{sarnak23letter} also sketch $\mathrm{GL}(1)$
  murmuration results using Poisson summation).

  \item Zubrilina~\cite{zubrilina23murmur}
  applies a trace formula of Eichler--Selberg  
  type due to Yamauchi~\cite{yamauchi73trace} (with corrections by Skoruppa and
  Zagier~\cite{skoruppazagier88trace}) to study murmurations of fixed weight,
  varying level holomorphic modular cuspforms.
  This was the first proved murmuration for $L$-functions not of
  $\mathrm{GL}(1)$ type, influencing several later works.

  \item Wang~\cite{wang2025murmurations} used similar analysis to Zubrilina, as
  well as character orthogonality and Poisson summation, to study murmurations
  of Hecke $L$-functions of quadratic imaginary number fields.

  \item Martin (mostly in~\cite{martin23rootnumber}, but
  see~\cite{kimballmartin2025murmur} for more murmuration discussion) also uses
  Yamauchi--Skoruppa--Zagier  
  to prove murmuration-like phenomena for local root numbers of modular forms.

  \item Booker, Lee, Seymour-Howell, Zubrilina, and the author~\cite{BLLDSHZ24}
  use an explicit version of the Selberg trace formula (coming from unpublished
  work of Strömbergsson~\cite{andreaspre} and appearing also
  in~\cite{andreithesis}) to prove murmurations of Maass forms of level $1$,
  weight $0$, and varying eigenvalue.
\end{enumerate}

Although the details differ, the broad mechanism behind many of these proofs is
similar.
With radical oversimplification: start with a trace formula, and then perform
sufficient averaging in each non-fixed aspect.

\begin{remark}
  Cowan~\cite{cowan2024murmurations} proposes a distinct method to prove
  murmurations assuming the Ratios Conjectures~\cite{conrey08ratios}.
  Instead of a trace formula, Cowan uses an explicit formula to relate sums of
  $\lambda(p)$ values to one-level density of the $L$-functions, which can be
  understood via logarithmic derivatives of the $L$-functions.
  He then applies the Ratios Conjectures to study these distributions.
\end{remark}

\begin{remark}
  Very recent work of Sawin and Sutherland~\cite{sawin2025murmurations} use
  Voronoi summation to study murmurations of elliptic curves ordered by naive
  height.
  This approach is also distinct from automorphic trace formulas and merits
  further investigation.
\end{remark}

\subsection{Flaws to trace formula approaches to murmurations}

There are problems when relying on trace formulas to study murmuration
phenomena.
Most notably, not enough trace formulas are known to explain all the
murmuration phenomena that can be experimentally observed (see Sutherland's
talk~\cite{sutherland24SCGPmurmuration} for an idea of the scope of
experimental observation).
The approaches in~\cite{ils2000lowlying, BBLLD23,
murmurDirichlet2023, zubrilina23murmur, wang2025murmurations,
martin23rootnumber, BLLDSHZ24} using trace formulas all fundamentally apply to
certain subsets of automorphic forms constructed around certain local
conditions.

No such approach will apply to sparse subsets of automorphic forms, such as
elliptic curves or genus $2$ curves.

But there are also qualitative problems.
For example, it seems very difficult to explain the similarities between the
arithmetically normalized murmuration~\eqref{eq:ils_murmuration} and the
unnormalized murmuration density~\eqref{eq:nu}, although both describe
correlation behavior for level $1$ holomorphic cusp forms of increasing weight.
Further, the murmuration density function for Maass forms of increasing
eigenvalue in~\cite{BLLDSHZ24} is \emph{exactly equal} to the density
function~\eqref{eq:nu} from~\cite{BBLLD23}, despite the fact that they come
from different families of $L$-functions.

\section{Explainable murmuration phenomena}

Practice shows that using a \emph{reasonable} trace formula and performing
\emph{reasonable} averages will \emph{probably} prove some murmuration
phenomena.
The structure of trace formulas is too broad to hope to make this concrete.

We conclude by describing several new murmurations that can \emph{probably} be
proven.

\subsection{Arithmetically weighted Maass forms}\label{ssec:arith_maass}

The Kuznetsov trace formula for level $1$ and weight $0$ relates coefficients
$\lambda_j(n)$ of Maass forms $\mu_j$ with eigenvalues $\frac{1}{4} + r_j^2$
and takes the form~\cite[Theorem~16.3]{iwanieckowalski04}
\begin{equation}
  \sum_{j \geq 1}
  \frac{\lambda_j(m) \lambda_j(n)}{\lVert \mu_j \rVert^2}
  \frac{h(r_j)}{\cosh \pi r_j}
  +
  \widehat{h}(0)
  =
  \sum_{c \geq 1}
  \frac{S(m, n; c)}{c}
  k^*(4 \pi \sqrt{mn} / c)
  +
  \delta_{[m = n]} k(0).
\end{equation}
In this expression, $h(\cdot)$ is a ``nice'' test function and all three of
$\widehat{h}(\cdot)$, $k(\cdot)$, and $k^*(\cdot)$ are weighted Fourier
transforms of $h(\cdot)$.

This is analogous to the Petersson Trace formula~\eqref{eq:petersson}.
It has the similar effect of having a relatively straightforward analysis in
terms of Kloosterman sums instead of class numbers, but unfortunately
introduces arithmetic weights $\lVert \mu_j \rVert^2 \asymp L(1, \mathrm{Sym}^2
\mu_j)$.

Using this with $m = 1, n = p$ should yield an arithmetically weighted
murmuration that is qualitatively similar to~\cite{BLLDSHZ24}.

This has an analogous relationship to~\cite{BLLDSHZ24}
as the murmuration result~\eqref{eq:ils_murmuration} from~\cite{ils2000lowlying}
has to~\cite{BBLLD23}; it's the same family of $L$-functions, but one has an
arithmetic normalization and one doesn't.
The murmuration densities in the nonnormalized cases of~\cite{BBLLD23}
and~\cite{BLLDSHZ24} are identical.
It is not clear whether one should expect strong similarity between the
arithmetic normalizations.

\subsection{Extending known murmurations}

The murmurations in~\cite{BBLLD23} and~\cite{BLLDSHZ24} apply to level $1$
holomorphic cuspforms and level $1$, weight $0$ Maass forms, respectively.
We should expect that both of these results can be generalized to general
level.

Both the Eichler--Selberg  
and Petersson trace formulas for holomorphic cuspforms are sufficiently
concrete in general.
Some of the work towards applying the Petersson trace formula to level $N$
forms appears in~\cite{ils2000lowlying}; it's almost certainly possible to
generalize that approach to prove the level $N$ version of the arithmetically
normalized murmuration in~\eqref{eq:ils_murmuration}.
Similarly, the trace formula in~\cite{child22} is explicit, and the techniques
from~\cite{BBLLD23} would generalize to arbitrary level.

In both cases, the resulting murmuration densities are essentially the same,
except for some local factors coming from primes dividing the level.

On the other hand, there \emph{are} major obstructions preventing
generalization of~\cite{BLLDSHZ24} to either general level or general weight
Maass forms.
The Selberg-Strömbergsson trace formula is only written down explicitly when
the level is squarefree.

Using the Kuznetsov trace formula (cf.\ \S\ref{ssec:arith_maass}), it should be
possible to prove arithmetically normalized murmurations for $L$-functions of
Maass forms on general level.
Further, one can establish an explicit Kuznetsov trace formula for any weight
$k$ (see for example~\cite[\S5]{dfi}), allowing one to consider murmurations
on arbitrary level and arbitrary weight as the eigenvalue aspect varies.

\begin{remark}
  Forthcoming work of Bober, Booker, Knightly, Lee, Seymour-Howell, and the
  author intends to give a sufficiently explicit trace formula to allow
  computing Maass forms.
  It should also enable proving more general murmurations, including arbitrary
  level and weight $1$ Maass forms.
\end{remark}

\subsection{Symmetric square murmurations}

Thus far, we have murmurations for $\mathrm{GL}(1)$ families and
$\mathrm{GL}(2)$ families.
It should be possible to investigate murmurations of symmetric square lifts of
holomorphic modular forms.
Recall that if $f$ is a holomorphic modular form with coefficients
$\lambda_f(n)$, then $\mathrm{Sym}^2(f)$ is a modular form on $\mathrm{GL}(3)$
with coefficients $\lambda_{\mathrm{Sym}^2(f)}(p) = \lambda_f(p^2)$.

Thus one can apply a $\mathrm{GL}(2)$ trace formula (such as
Eichler--Selberg  
or Petersson) to traces of coefficients $\lambda_f(p^2)$ to study murmuration
behavior across symmetric square lifts.

Analogously, one can try to use the Kuznetsov trace formula to study
murmurations in symmetric square lifts of Maass forms by studying correlations
in $\lambda_j(p^2)$.

\begin{remark}
  The root numbers of symmetric square lifts of $\mathrm{SL}(2, \mathbb{Z})$
  modular forms all have root number $1$.
  Thus it doesn't make sense to examine ``correlation'' with root number.
  Nonetheless, murmuration phenomena of the shape~\eqref{eq:murmuration_shape}
  may still exist.
\end{remark}

\begin{remark}
  More generally, one could consider arbitrary sparse subsequences.
  Steven Creech is studying statistical properties of $\lambda_f(q(n))$, where
  $q(n)$ is an irreducible quadratic polynomial. It would be desirable to study
  $\lambda_f(p^m)$ to understand $\mathrm{Sym}^m(f)$ and thus higher
  $\mathrm{GL}(n)$ behavior, but initial heuristic estimates suggest this would
  be very challenging to carry out in practice.
\end{remark}

\subsection{Murmurations for families without the Riemann Hypothesis}

Finally, we outline an approach to prove murmurations for a family of Dirichlet
series that don't satisfy their Riemann hypothesis.
To a weight $k + \tfrac{1}{2}$ cuspidal modular form $f$ (as
in~\cite{shimura}), one can associate a Dirichlet series $D(s, f)$ that has a
functional equation
\begin{equation}
  N^s D(s, f) G(s) = \varepsilon N^{1-s} D(1-s, \widetilde{f}) G(1-s)
\end{equation}
(cf.\ the standard functional equation~\eqref{eq:feq}), analytic continuation,
and other properties of standard automorphic $L$-functions --- except that they
don't have Euler products and have nontrivial zeros outside the critical line.

Nonetheless, a version of the Petersson trace formula applies (see for example
the appendix to~\cite{KLDWS23} for a related application).
Much of the analysis seems similar to~\cite{BBLLD23}.
Initial sketches suggest that GRH for Dirichlet $L$-functions and using a
smooth cutoff function $\Phi(\cdot)$ instead of a sharp cutoff function should
be sufficient.

\begin{remark}
  The author intends to investigate this in later work with the authors
  of~\cite{BBLLD23, BLLDSHZ24}.
\end{remark}

\bibliographystyle{alpha}
\bibliography{bibfile}

\end{document}